\documentstyle[12pt]{article}

\def\ra{\rightarrow}

\def\ss{\subseteq}

\def\e{\epsilon}

\def\O{\Omega}

 \def\HollowBox #1#2{{\dimen0=#1 \advance\dimen0 by -#2       
       \dimen1=#1 \advance\dimen1 by #2                       
        \vrule height #1 depth #2 width #2                    
        \vrule height 0pt depth #2 width #1                   
        \llap{\vrule height #1 depth -\dimen0 width \dimen1}%
       \hskip -#2                                             
       \vrule height #1 depth #2 width #2}}                   
              
\font\teneufm=eufm10
\font\seveneufm=eufm7
\font\fiveeufm=eufm5
\newfam\eufmfam
\textfont\eufmfam=\teneufm
\scriptfont\eufmfam=\seveneufm
\scriptscriptfont\eufmfam=\fiveeufm

\newfam\msbfam
\font\tenmsb=msbm10 scaled \magstep1  \textfont\msbfam=\tenmsb
\font\sevenmsb=msbm7 scaled \magstep1  \scriptfont\msbfam=\sevenmsb
\font\fivemsb=msbm5  scaled \magstep1   \scriptscriptfont\msbfam=\fivemsb
\def\Bbb{\fam\msbfam \tenmsb}

\def\RR{{\Bbb R}}
\def\CC{{\Bbb C}}

\newtheorem{theorem}{Theorem}

\begin{document}

\begin{center}
\Large \bf Calculation and Estimation of the
Poisson Kernel 
\bigskip  \\
\small \rm by 
\medskip \\
\normalsize \rm Steven G. Krantz\footnote{Author
supported in part by NSF Grant DMS-9988854.}

\end{center}
\vspace*{.15in}

\begin{quote}
{\small \bf Abstract:}  \sl We provide a simple
method for obtain boundary asymptotics of the
Poisson kernel on a domain in $\RR^N$.
\end{quote}

\setcounter{section}{-1}

\section{Introduction}

Let $\Omega \subseteq \RR^N$ be a connected open set---called
a {\it domain}.  It is a matter of considerable interest
to estimate the size of the Poisson kernel 
$P_\Omega(x,t) = P(x,t)$ of $\Omega$.
Here, and throughout this paper, $x \in \Omega$
and $t \in \partial \Omega$.

In case $\Omega$ has a large group of symmetries, then
it is often possible to calculate $P_\Omega$
explicitly.  For example,
\begin{itemize}
\item The Poisson kernel of the disc $D \ss \RR^2$ is
$$
P_D(x,t) = \frac{1}{2\pi} \cdot \frac{1 - |x|^2}{|x - t|^2} \, .
$$
\item The Poisson kernel for the upper halfplane
$$
U^2 = \{(x_1, x_2) \in \RR^2: x_2 > 0\}
$$
is given by
$$
P_{U^2}(x,t) = \frac{1}{\pi} \cdot \frac{x_2}{(x_1 - t)^2 + x_2^2} \, .
$$
\item The Poisson kernel for the unit ball $B \ss \RR^N$
is given by
$$
P_B(x,t) = \frac{\Gamma(N/2)}{2 \pi^{N/2}} \cdot \frac{1 - |x|^2}{|x - t|^N} \, .
$$
Here $\Gamma$ is the classical gamma function.
\item The Poisson kernel for the upper halfspace $U^{N+1} \equiv
\{x = (x_1,\dots, x_{N+1}) \in \RR^{N+1}: x_{N+1} > 0\}$
(with $x = (x_1,\dots, x_{N+1}) = (x', x_{N+1})$) is given by
$$
P_{U^{N+1}}(x,t) = c_N \frac{x_{N+1}}{([x' - t]^2 + x_{N+1}^2)^{[N+1]/2}} \,
$$
where
$$
c_N = \frac{\Gamma([N+1]/2)}{\pi^{[N+1]/2}} \, .
$$
\end{itemize}

For purposes of studying the Schauder estimates for
the Dirichlet problem, for studying the (nontangential)
boundary behavior of harmonic functions, and for studying
potential theory, one needs to have size estimates for
the Poisson kernal on a fairly general domain (say
a bounded domain with $C^2$ boundary).

The standard asymptotic is
$$
P_\Omega(x,y) \approx \frac{\delta(x)}{|x - y|^N} \, .	\eqno (*)
$$
Here $\delta(x) \equiv \delta_\Omega(x)$ is the distance
from $x \in \Omega$ to $\partial \Omega$.  This estimate, together
with analogous estimates for the derivatives of $P_\Omega$,
suffices for most applications.  It is our purpose in
this paper to give an efficient and elementary method for
proving $(*)$.  At the end of the paper we shall also
sketch an argument for obtaining the cognate estimate
for derivatives of the Poisson kernel.

There are a number of methods for deriving estimates as
we have described, though none of them is well known.
After all, the harmonic analysis of domains in space
is a fairly new field, and many of the techniques
are only recently born.  Classical studies, in dimension
two only, appear in [KEL].  In the reference [KRA1],
we present an argument based on Kelvin reflection
of harmonic functions and comparisons by way of
the maximum principle.  These arguments were developed
by Norberto Kerzman (personal communication).
They were presented with Kerzman's permission.  They
are intricate, and we shall not repeat them here.

Another natural method for developing an asymptotic
expansion for the Poisson kernel is to use Fourier
integral operators.  To wit, let us suppose for simplicity
that $\Omega$ is topologically trivial and has smooth
boundary.  Let $\Phi:\overline{\Omega} \ra \overline{B}$ 
be a diffeomorphism of the closure of $\Omega$ with the
closure of $B$.  Then one can compare the true Poisson
kernel on $\Omega$ with the pullback of the Poisson
kernel from $B$ under the mapping $\Phi$.  The result
is the required asymptotic expansion (see [BMS], where
a similar technique is used to obtain an asymptotic
expansion for the Bergman kernel and the Szeg\"{o} kernel).

In the present paper we use a method that has come to
be known as ``scaling'' to produce the estimates $(*)$ for
the Poisson kernel.  This is a methodology that has been
developed extensively in the study of automorphism groups
of domains in $\CC^n$---see [ISK].  It has also been used in harmonic
analysis to obtain information about reproducing kernels
(see [NRSW], where it was used to study the Szeg\"{o} kernel).
The advantages of this approach are that {\bf (i)} it
is quite elementary and straightforward and {\bf (ii)} it
can be applied to a variety of reproducing kernels in
many different circumstances.  Thus the techniques
presented here should find utility in a number of different
contexts.

It is a pleasure to thank Kang-Tae Kim, Richard Rochberg,
and Norm Levenberg for helpful conversations.  Kim has taught
me much of what I know about scaling.  Coifman and Rochberg
[COR] prove estimates much like the ones presented here,
but on the ball for a Bergman space with weights.  Levenberg
and Yamaguchi [LEY] use a scaling method similar to the one here
to estimate a reproducing kernel from another context.

\section{The Main Result}

For the remainder of the paper, let $\Omega \ss \RR^N$ be
a bounded domain with $C^2$ boundary.  This means
that there is a $C^2$, real-valued function $\rho$ such
that
$$
\Omega = \{x \in \RR^N: \rho(x) < 0\}
$$
and $\nabla \rho \ne 0$ on $\partial \Omega$.  Thus
$\partial \Omega$ is a regularly imbedded $C^2$ hypersurface
in $\RR^N$.  

\begin{theorem} \sl
Let $\Omega \ss \RR^N$ be a bounded domain with $C^2$ boundary.
Let $P: \Omega \times \partial \Omega \ra \RR^+$ be the 
Poisson kernel for $\Omega$.  Then there are constants
$c_1, c_2 > 0$ such that
$$
c_1 \cdot \frac{\delta(x)}{|x - y|^N} \leq P(x,y) \leq
      c_2 \cdot \frac{\delta(x)}{|x - y|^N} \, .    \eqno (\star)
$$
\end{theorem}

The remainder of this section is devoted to the proof
of this theorem.  In the last section of the paper we shall remark
on how to obtain a similar asymptotic for the
derivatives of $P$.  For convenience, we write
$$
P(x,y) \approx \frac{\delta(x)}{|x - y|^N} 
$$
instead of $(\star)$. 

Notice before we begin that, if $K$ is a compact set in $\Omega$,
then the estimate we seek is trivial for $x \in K$ and $y \in \partial
\O$.  For then $|x - y| \geq c > 0$, $\delta(x)$ is bounded above,
and we get a universal bound above and below on $\delta(x)/|x -
y|^N$.  A similar comment applies if $x$ is near the boundary
and $y$ is far from $x$.  So we may concentrate our attention
on $x$ near the boundary and $y$ near $x$.

Now fix a point $P \in \partial \Omega$ and a point $P^0 \in \Omega$
such that the segment $\overline{P^0 P}$ is normal to
the boundary at $P$.  We shall dilate coordinates with
center $P^0$.  We assume that $P^0$ is close to $\partial \Omega$---within
a tubular neighborhood of the boundary---and we set
$\epsilon = \hbox{dist}(P^0,P)$.  
We assume that coordinates have been rotated and centered
so that
\begin{enumerate}
\item[{\bf (a)}]  The point $P$ is the origin $(0,0,\dots,0)$;
\item[{\bf (b)}]  The normal direction $\overrightarrow{P P^0}$
is the positive $x_N$-direction.  
\end{enumerate}

For a point $x \in \RR^N$, we write $x = (x_1,\dots, x_N)$.
We set $P^0 = (P^0_1, \dots, P^0_N)$.  With the normalization
of coordinates, in fact $P^0 = (P_1^0, P_2^0, \dots, P_N^0) 
= (0,\dots, 0, +\epsilon)$.
Now define
$$
\Phi_\epsilon(x) = \left ( \frac{x_1}{\e}, \frac{x_2}{\e},
                       \dots, \frac{x_N}{\e} \right ) \, .
		       $$
Observe particularly that the mapping $\Phi_\epsilon$ sends
the point $P^0$ to $(1,0,\dots, 0)$.  

The first thing to notice is that, in a natural
sense,
$$
\lim_{\e \ra 0^+} \Phi_\e(\Omega) = U^N \, .
$$
To see this, we first check that if the defining function $\rho$,
expanded about the point $P$, is given by
$$
\rho(x) = \sum_{j=1}^N a_j^1 x_j + \sum_{j,k=1}^N a_{jk}^2 x_j x_k 
   + \cdots = - x_N + \sum_{j,k=1}^N a_{jk}^2 x_j x_k \, .
$$
(of course note that $\rho(P) = 0$) then
$$
\rho_\e(s) \equiv \frac{1}{\epsilon} \cdot \left [
    \rho \circ \Phi_\epsilon^{-1} (s) \right ] = \frac{1}{\epsilon} \cdot
				\left [ - \e s_N + \sum_{j,k = 1}^N a_{jk}^2 \epsilon^2 s_j s_k + \cdots \right ]
   = - s_N + \epsilon \cdot \left [ 
         \sum_{j,k = 1}^N a_{jk}^2 s_j s_k + \cdots \right ] \, .		       
$$     
Plainly, as $\e \ra 0$, the transferred defining
function $\rho_\e$ tends to the linear defining function 
$\rho_0(s) \equiv - s_N$.
In other words, the domains $\Phi_\e(\Omega) \equiv \Omega_\e$ converge
(in an appropriate sense) to the standard halfspace.
This last information is useful because we know the
Poisson kernel for a halfspace.

Now we may take advantage of the facts accrued by 
setting $\O_\e = \Phi_\e (\O)$, letting $d\sigma$ be $(N-1)$-dimensional
area measure on $\partial \O$, 
$d\sigma_\e$ to be $(N-1)$-dimensional area measure on $\partial \O_\e$,
and taking $f$ to be a
function that is continuous on $\overline{\O}_\e$ and
harmonic on $\O_\e$.
Further, we let $x \in \O$ and set $s = \Phi_\e(x)$.
Then we calculate that
\begin{eqnarray*}
f(s) & = & f(\Phi_\e(x)) \\
     & = & \int_{\partial \O_\e} P_{\O_\e} (\Phi_\e(x), t) f(t) \, d\sigma_\e(t)  \\
     & = & \int_{\partial \O} P_{\O_\e} (\Phi_\e(x), \Phi_\e(\tau))
              f(\Phi_\e(\tau)) \, \hbox{det} \, \hbox{Jac} \Phi_\e(\tau)
	            \, d\sigma(\tau) \, .
\end{eqnarray*}
It is crucial to note here that the integral is over an 
$(N-1)$-dimensional hypersurface, and hence the
Jacobian determinant is that of an $(N-1) \times (N-1)$ matrix.

Now let us write
$$
K_\e(x, \tau) = P_{\O_\e}(\Phi_\e(x), \Phi_\e(\tau)) \cdot
		    \hbox{det} \, \hbox{Jac} \Phi_\e(\tau) 
	      = \e^{-(N-1)} \cdot P_{\O_\e}(\Phi_\e(x), \Phi_\e(\tau)) \, .
$$
We thus have the equation
$$
f \circ \Phi_\e(x) = \int_{\partial \O} P_\O(x,\tau) 
    [ f\circ \Phi_\e(\tau)] \, d\sigma(\tau)
	   = \int_{\partial \O} K_\e(x,\tau) [f \circ \Phi_\e(\tau)] \, d\sigma(\tau) \, .
  \eqno (\dagger)
$$
Since this identity holds true for any choice of continuous
$f$ on the boundary of $\Omega_\e$ (with unique harmonic extension
to $\Omega_\e$), we may conclude that
$$
P_\O(x,\tau) = K_\e(x,\tau) \, .  \eqno (\ddagger)
$$

The identity $(\star)$ is the key to our result, for
we know asymptotically what $K_\e$ looks like.
In particular, we know (see [KRA1, Section 1.3]) 
on any smoothly bounded domain
$U$ that the Poisson kernel is a normal derivative
of the Green's function:
$$
P_U(x,y) = \frac{\partial}{\partial \nu_y} G_U(x,y) \, .
$$
And the Green's function, in turn, is the solution on
$U$ of the Dirichlet problem with boundary data
the Newton potential $\Gamma_N( \ \cdot \ - x)$.

Now with $P, P^0$ fixed as before, let $W$ be a small,
smoothly bounded, topologically trivial domain  
with these properties:
\begin{enumerate}
\item[{\bf (a)}]  $\displaystyle W \subseteq \Omega$;
\item[{\bf (b)}]  $\displaystyle P^0 \in W, P \in \partial W$;
\item[{\bf (c)}]  $\partial W \cap \partial \Omega$ is a 
relative neighborhood of $P$ in $\partial \Omega$.
\end{enumerate}
Easy Schauder estimates, and the discussion in the
preceding paragraph, show that we may obtain our estimate
$(\star)$ by studying the cognate question on $W$ (details of
this type of argument may be found in [APF]).

Now the key observation at this point is that, when $\e > 0$
is small, then the Poisson
kernel for $\Phi_\e(W)$ at interior points of the line segment
$\Phi_\e(\overline{P P^0})$ is very near to the Poisson
kernel of the upper half space $U^N$ at those same points.
The reason, of course, is that if $\rho_1$ is the defining
function for $U^N$ and $\rho_2$ is the defining
function for $\Phi_\e(W)$ then there is a diffeomorphism
$\lambda$ so that $\rho_1 = \rho_2 \circ \lambda$ near 0
and $\|\lambda - \hbox{id}\|_{C^1}$ is small.  Referring
again to the construction of the Poisson kernel above, the
claim follows.  

As a result, we may calculate the Poisson kernel on $\Omega$
by instead calculating the kernel on $W$.  In turn, it
then suffices to calculate the kernel on $U^N$.  
Thus we see that, for $x$ on the interior of the line segment $\overline{P P^0}$,
\begin{eqnarray*}
K_\epsilon(x, \tau) & = &\e^{-(N-1)} \cdot P_{\O_\e}(\Phi_\e(x), \Phi_\e(\tau))  \\
             & \approx & \epsilon^{-(N-1)} \frac{\Phi_\e(x)_N}{(|\Phi'_\e(x) 
                - \Phi_\e(\tau)|^2 + [\Phi_\e(x)_N]^2)^{N/2}}  \\
             & = & \e^{-(N-1)} \cdot \frac{x_N/\e}{(|x'/\e - \tau/\e|^2
                      + [x_N/\e]^2)^{N/2}}  \\
             & = & \frac{x_N}{(|x' - \tau|^2 + [x_N]^2 )^{N/2}} \, .
\end{eqnarray*}
Unraveling the notation, we find that we have proved the
approximation $(\star)$.

\section{Estimates for the Derivatives of the Poisson Kernel}

The argument to obtain an asymptotic for a derivative
of the Poisson kernel is nearly the same as the ones used above.  The result
we seek is
$$
\nabla_x^k P(x,y) \approx \frac{\delta(x)}{|x - y|^{N+k}} \, . \eqno (**)
$$
The crux of the argument is the analog of equation $(\dagger)$.
For the present application, that equation now becomes
$$
\nabla_x^k f \circ \Phi_\e(x) = \nabla_x^k \int_{\partial \O} P_\O(x,\tau) 
    [ f\circ \Phi_\e(\tau)] \, d\sigma(\tau)
	   = \int_{\partial \O} \nabla_x^k K_\e(x,\tau) 
          [f \circ \Phi_\e(\tau)] \, d\sigma(\tau) \, .
$$

It follows as before that we have the identity
$$
\nabla_x^k P_\Omega(x, \tau) = \nabla^k_x K_\e(x, \tau) + {\cal E} \, ,
$$
where ${\cal E}$ is an error term that is a polynomial (and hence
is of no interest for our estimates).  The remainder of
the derivation of $(**)$ is as before.
\vspace*{.2in}

\leftline {\Large \sc References}

\begin{enumerate}

\item[{\bf [APF]}]  L. Apfel, Localization properties and
boundary behavior of the Bergman kernel, thesis, 
Washington University in St.\ Louis, 2003.

\item[{\bf [BMS]}]  L. Boutet de Monvel and J. Sj\"{o}strand, Sur la
singularit\'{e} des noyaux de Bergman et Szeg\"{o}, {\em Soc.\ Mat.\ de
France Asterisque} 34-35(1976), 123-164.

\item[{\bf [COR]}]  R. Coifman and R. Rochberg, Representation
Theorems for Holomorphic and Harmonic Functions in $L^p$.
{\it Representation Theorems for Hardy Spaces}, pp. 11--66,
Asterisque 77, Soci\'{e}t\'{e} Math\'{e}matique de France,
Paris, 1980. 

\item[{\bf [ISK]}]  A. Isaev and S. G. Krantz, 
Domains with non-compact automorphism group:  a survey,
{\it Advances in Math.}\ 146(1999), 1--38.

\item[{\bf [KEL]}]  Keldych Lavrentieff,
Sur une evaluation de la fonction de 
Green, {\em Doklady Acad.\ USSR} 24(1939), 102.

\item[{\bf [KRA1]}]  S. G. Krantz, {\it Function Theory of
Several Complex Variables}, $2^{\rm nd}$ ed., American
Mathematical Society, Providence, RI, 2001. 

\item[{\bf [LEY]}]  N. Levenberg and H. Yamaguchi,
The metric induced by the Robin function, {\it Mem.\ Amer.\ Math.\ Soc.}\
92(1991), viii$+$156.

\item[{\bf [NRSW]}]  A. Nagel, J. P. Rosay, E. M. Stein, and S. Wainger,
Estimates for the Bergman and Szeg\"{o} kernels in $\CC^2$, 
{\em Ann.\ Math.}\ 129(1989), 113-149. 

\end{enumerate}

\end{document}